\documentclass{IEEEtran4PSCC}

\usepackage{cite}
\usepackage{enumitem}
\usepackage{pdflscape}
\usepackage{comment}
\usepackage[hidelinks]{hyperref}
\usepackage{caption}
\ifCLASSINFOpdf
   \usepackage[pdftex]{graphicx}
 
\else

   \usepackage[dvips]{graphicx}
  
\fi

\usepackage[cmex10]{amsmath}

\usepackage{array}

\usepackage{url}

\usepackage{longtable}
\usepackage{booktabs}
\usepackage{graphicx} 
\usepackage{amsmath}
\usepackage{amssymb}
\usepackage{multirow}
\usepackage{siunitx}
\usepackage{caption}
\usepackage{float}

\newcommand{\myred}[1]{#1}
\renewcommand{\myred}[1]{{\color{red}#1}} 
\newcommand{\myblue}[1]{{\color{blue}#1}}

\newcommand{\answer}[1]{{\color{black}#1}}

\usepackage{geometry}
\makeatletter
\let\old@ps@headings\ps@headings
\let\old@ps@IEEEtitlepagestyle\ps@IEEEtitlepagestyle
\def\psccfooter#1{%
    \def\ps@headings{%
        \old@ps@headings%
        \def\@oddfoot{\strut\hfill#1\hfill\strut}%
        \def\@evenfoot{\strut\hfill#1\hfill\strut}%
    }%
    \def\ps@IEEEtitlepagestyle{%
        \old@ps@IEEEtitlepagestyle%
        \def\@oddfoot{\strut\hfill#1\hfill\strut}%
        \def\@evenfoot{\strut\hfill#1\hfill\strut}%
    }%
    \ps@headings%
}
\makeatother

\begin{document}

\title{Integrated Investment and Operational Planning for Sugarcane-Based Biofuels and Bioelectricity under Market Uncertainty}

\author{\IEEEauthorblockN{Carolina Monteiro\IEEEauthorrefmark{1}\IEEEauthorrefmark{2},
Bruno Fanzeres\IEEEauthorrefmark{2},
Rafael Kelman\IEEEauthorrefmark{1}, 
Raphael Araujo Sampaio\IEEEauthorrefmark{1}, \\ 
Luana Gaspar\IEEEauthorrefmark{1},
Lucas Bacellar\IEEEauthorrefmark{1}, and
Joaquim Dias Garcia\IEEEauthorrefmark{1},}
\IEEEauthorblockA{\IEEEauthorrefmark{1} PSR\\
 \{carolina, rafael, rsampaio, luana, lucasbacellar, joaquim\}@psr-inc.com}
\IEEEauthorblockA{\IEEEauthorrefmark{2} Industrial Engineering Department, PUC-Rio\\ bruno.santos@puc-rio.br}
}

\maketitle

\begin{abstract}
Sugarcane biomass is a strategic resource for the energy transition, particularly in Brazil, where it underpins electricity and ethanol production. Investment planning is challenged by diverse production pathways, price volatility, and feedstock variability. This work develops a two-stage stochastic optimization model integrating investment and operational decisions for sugarcane facilities. The model aims to support robust planning for diversified biomass plants, aiding the sector's decarbonization. The first stage defines capacity expansion under economies of scale through a power-law cost function. The second stage defines operational decisions under price and feedstock uncertainty, modeled via scenarios and Conditional Value-at-Risk. \answer{From an investor's perspective, the objective is to minimize risk-adjusted net costs. In addition to its methodological contributions, this work also provides an open-source implementation of the proposed capacity expansion planning framework, referred to as \textit{OptBio}.} A Brazilian case study shows risk-neutral strategies favor sugar/ethanol but are vulnerable, whereas risk-averse strategies promote diversification. Sensitivity analyses indicate biomethane and hydrogen may become viable with favorable prices, while biochar boost productivity and profitability.

\end{abstract}


\thanksto{\noindent This study was partially supported by the Coordenação de Aperfeiçoamento de Pessoal de Nível Superior – Brasil (CAPES) – Finance Code 001. Bruno Fanzeres gratefully acknowledges the financial support from CNPq, projects 309064/2021-0 and 302272/2025-9, and FAPERJ, project E-26/204.553/2024.}

\section{Introduction}

    Addressing the escalating climate crisis is a global priority that requires coordinated efforts across multiple sectors (e.g., energy, industry, transport, and land use) to reduce greenhouse gas (GHG) emissions, prevent dangerous levels of global warming, and support the transition to a net-zero, low-carbon economy. Within the energy sector, expanding the use of low-carbon renewable sources for diverse energy production is widely recognized as a key and feasible strategy for achieving sustainable goals and supporting the desired energy transition. Among these renewable sources, biomass has emerged as a strategic resource, offering a reliable alternative to carbon-intensive fossil fuels. In Brazil, for example, sugarcane biomass is currently the main renewable source in the total energy supply, and, along with other biomass sources, accounted for over 60,000 gigawatt-hours of electricity generation in 2024 \cite{EPE2025}, of which 24,417 gigawatt-hours was exported to the grid \cite{ons_energy_generation}. Furthermore, in the transport sector, ethanol (primarily derived from sugarcane) and biodiesel together account for approximately 25\% of total energy consumption \cite{EPE2025}. In fact, recent technological advancements have further expanded this potential, enabling the production of biogas, biomethane, and advanced biofuels \cite{borjesson2008, Wang2016, Aguiar2025}, which are key to decarbonizing hard-to-electrify sectors such as aviation and heavy-duty transport.

    However, the variety of operational configurations and the multiple derivatives/products that can be produced from biomass (particularly from sugarcane) create significant planning and management challenges for investors in this sustainable energy source. These challenges include assessing an adequate investment plan for the production facility, accurately modeling the complexity of the production chain, and efficiently managing the uncertainty associated with key derivative prices (e.g., sugar and ethanol) and feedstock availability \cite{Goldenberg2008, Henniges2004}. More specifically, appropriate investment decisions should consider diverse production pathways and the uncertainties related to both price fluctuations and feedstock availability. They should also be driven by flexible operational models that integrate agricultural, industrial, and market variables (e.g., the flexibility to switch between ethanol and sugar depending on relative profitability) \cite{Goldenberg2008}. Nevertheless, due to these challenges, the technical and scientific literature often simplifies methodologies for devising investment plans in sugarcane biomass facilities. This simplification typically narrows the range of considered products, inadequately represents uncertainty factors, oversimplifies operational aspects of the production chain, or decouples investment decisions from plant operations.

    \answer{Therefore, this work proposes a capacity expansion planning model for sugarcane biomass production facilities that integrates investment and operational decisions while explicitly accounting for uncertainties in both prices and product availability, and accurately representing the complexity of the production chain. Structurally, given the nature of these uncertainties, the proposed decision-making framework is framed as a two-stage stochastic program, in which the first stage determines the investment plan and the second stage optimizes the corresponding operational decisions.} Uncertainty is represented by a set of scenarios for both the prices of each derivative and the availability of products within the production chain for a typical year of operation. To account for the investor's risk profile, we employ the quantile-based risk measure known as Conditional Value-at-Risk (CVaR) \cite{Rockafellar2002_CVaRforGeneralLossDist}, leveraging its intuitive specification and interpretation, as well as its high computational tractability. Conceptually, the proposed methodology aims to balance investment costs with risk-adjusted operational expenditures, the latter being captured through a convex combination of the expected value and the CVaR of the operating net costs across the scenarios of price and availability. \answer{In the first stage, the model determines the capacity and associated investment for each candidate plant, accounting for economies of scale through the capacity-cost scaling power-law relationship \cite{chilton, tribe}. The second stage then defines the plant operations for a representative (typical) year\footnote{\answer{In this work, we assume that the complete investment plan is determined at the beginning of the operational horizon, characterized by a representative (typical) year. Consequently, aspects such as multi-year investment policies, the dynamic evolution of uncertainties, and asset degradation and depreciation are not explicitly modeled. Furthermore, we restrict our analysis to the design of capacity expansion plans that explicitly account for facility operations. Accordingly, the annual sugarcane harvest is treated as an exogenous, fixed condition. As a result, aspects related to the agricultural side of the system (e.g., land-use decisions, crop-switching, and the detailed logistics and costs of straw collection) are not explicitly modeled.}} \cite{Tanaka2022_EconomicsofPwrSyst, Ribeiro2024_RepresentativeDays_EnergyAI} based on the installed capacities and the realization of price and availability scenarios, with the objective of minimizing total net costs, composed of operational expenditures offset by sales revenue. We argue that this level of temporal aggregation (namely, the use of a representative (typical) year) is sufficiently accurate for informing and determining investment plans for sugarcane biomass production facilities, given the strong annual cyclicality of sugarcane biomass–related activities. Furthermore, for expositional clarity and to better leverage available information and data, and recognizing the relative stability of operations throughout the year, we represent within-year operations using a single time period.}

    \answer{The numerical experiment demonstrates the model's applicability. The results indicate that the widely adopted practice in Brazil of investing in biomass-to-electricity capacity yields stable returns but remains exposed to fluctuations in sugar and ethanol prices. In contrast, incorporating risk aversion underscores the benefits of diversification through technologies such as Alcohol-to-Jet (AtJ), which improves the overall risk–return trade-off. Additionally, waste-to-energy pathways (e.g., biomethane and hydrogen) and productivity-enhancing options (e.g., biochar) can generate new revenue streams, enhance profitability, and reinforce the role of biomass facilities in the transition toward low-carbon energy systems.}

\subsection{Objectives and contributions regarding the existing literature}

Over the past decades, the challenges associated with planning and operating biomass-based fuel production facilities have received significant attention in research. For example, \cite{mutran2020risk} proposes a risk-aware optimization model to guide investment decisions in the Brazilian sugarcane industry, focusing on standard sugar mill products (sugar, ethanol, fertilizer, and electricity) and evaluating the most financially suitable processes for producing these products. Additionally, \cite{gilani2021multi} presents a robust optimization model for designing a sugarcane-based ethanol supply chain, emphasizing strategic investment decisions (location and capacity) and operational decisions (production, transportation, and inventory) under uncertainty regarding demand and export prices. In \cite{jonker2016supply}, a deterministic optimization model is used to address the supply chain of first-generation (1G) ethanol from sugarcane and second-generation (2G) ethanol from eucalyptus, leveraging a Mixed-Integer Linear Programming (MILP) model to evaluate various expansion strategies for determining the most suitable locations and capacities for processing plants. Furthermore, \cite{woo2016optimization} introduces a deterministic model for optimizing the design and operation of a biomass-to-hydrogen supply chain, considering biomass importation, inventory strategies, and fluctuating hydrogen demand.

    In \cite{dutenkefer2018insertion}, the authors used robust optimization to explore adding biogas production to sugarcane mills. Their model, based on a fixed plant layout, allocated sugarcane optimally by using historical prices to account for uncertainty. Expanding on the model in \cite{mutran2020risk}, study \cite{watson2025case}  looked at the same core products but included a new process: ethanol upgrading, which converts ethanol into valuable fuels like gasoline, Sustainable Aviation Fuel (SAF), and diesel. Unlike \cite{mutran2020risk}, the study in \cite{watson2025case} focuses exclusively on scenario-dependent operational decisions, not investment choices. It allows the production mix to vary based on different price scenarios. The main goal is to determine how operational decisions shift in response to relative product prices, specifically identifying when SAF becomes the most profitable product. \answer{Table \ref{tab:literature_comparison} provides a structured summary of the representative literature, organized according to key modeling features.}

    \begin{table}[!ht]
        \centering
        \caption{\answer{Summary of the representative literature, organized according to key modeling features.}}
        \label{tab:literature_comparison}
        \renewcommand{\arraystretch}{1.30}
        \resizebox{\columnwidth}{!}{%
            \begin{tabular}{r ccc}
                \hline
                \multirow{2}{*}{ } & \textbf{Integrated Investment} & \textbf{Risk}         & \textbf{Economies} \\
                                   & \textbf{\& Operations}         & \textbf{Management}   & \textbf{of Scale} \\
                \hline
                Mutran et al. \cite{mutran2020risk} & $\checkmark$ & $\checkmark$ & $\times$ \\
                
                Gilani and Sahebi \cite{gilani2021multi} & $\checkmark$ & $\checkmark$ & $\times$ \\
                
                Jonker et al. \cite{jonker2016supply} & $\checkmark$ & $\times$ & $\times$ \\
                
                Woo et al. \cite{woo2016optimization} & $\checkmark$ & $\times$ & $\times$ \\
                
                Dutenkefer et al. \cite{dutenkefer2018insertion} & $\times$ & $\checkmark$ & $\times$ \\
                
                Watson et al. \cite{watson2025case} & $\times$ & $\times$ & $\times$ \\
                
                \hline
                \textbf{This work} & $\checkmark$ & $\checkmark$ & $\checkmark$ \\
                \hline
            \end{tabular}
        }
    \end{table}

Despite the relevance of the aforementioned studies, current works have largely concentrated on isolated challenges within the sugarcane biomass production chain and its capacity expansion. A comprehensive investment plan that captures key uncertainties and accurately reflects the complexity of the production chain remains underexplored. This work fills that gap by proposing a capacity expansion planning model for sugarcane biomass facilities that integrates both investment and operational decisions. The model explicitly accounts for the primary uncertainty factors impacting the business (i.e., derivative prices and product availability) and represents a wide range of products, plants, and processes, extending beyond the usual focus on basic outputs like sugar and ethanol. Additionally, the methodology incorporates nonlinear economies of scale that link plant capacity expansion to Capital Expenditures (CAPEX), enhancing the model’s precision and suitability for real-world applications.

To further enhance the accessibility and applicability of the proposed methodology, an open-source software package based on the optimization model, OptBio \cite{OptBio}, is developed. The software is implemented in the Julia language \cite{Julia} using the JuMP mathematical optimization framework \cite{JuMP}, which facilitates the use of both commercial and open-source solvers, such as HiGHS \cite{HiGHS}. Fully documented and user-friendly, the package allows researchers and professionals to customize their study cases through a programmatic interface, storing data in an SQLite database \cite{sqlite2020hipp}. This open-source solution ensures that the results presented in this work can be reproduced and adapted by others in the field, fostering broader collaboration and transparency.
Finally, to demonstrate the applicability of the proposed methodology and offer insights for the sugarcane sector, we present a numerical experiment based on the Brazilian context, incorporating multiple plants, processes, and products. The case study includes electricity, hydrogen, biofuels, and a diverse range of low-carbon sugarcane-derived products, extending the typical scope of applications found in the scientific and technical literature and providing a more comprehensive view of production alternatives.

\begin{figure*}[hbtp]
    \centering
    \includegraphics[width=\textwidth]{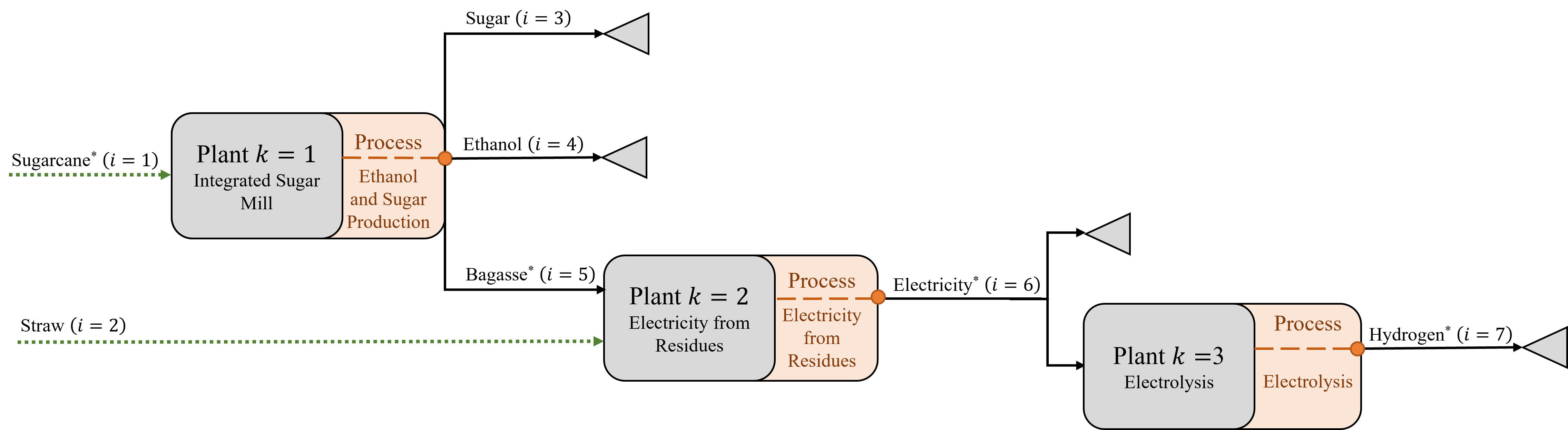}
    \caption{Production chain of an example facility composed of $|\mathcal{K}| = 3$ plants that involves a total of $|\mathcal{I}| = 7$ products.}
    \label{fig:Illust_Prod}
\end{figure*}

\section{Capacity Expansion Planning Model} \label{Sec::CapExpModel}

The primary objective of this work is to propose a capacity expansion planning model for sugarcane biomass production facilities that integrates both investment and operational decisions. The model explicitly accounts for uncertainties in market prices and product availability, while also capturing the complexity of the production chain and the nonlinear economies of scale associated with infrastructure investments. Structurally, the proposed model follows a two-stage decision-making framework \cite{Tanaka2022_EconomicsofPwrSyst}. In the first stage, an investment plan for production plants is determined under uncertainty in both market prices and product availability. In the second stage, given the planned capacity expansion, the model optimizes the operation of the sugarcane biomass production facility for a representative (typical) year \cite{Tanaka2022_EconomicsofPwrSyst}.

From an operational standpoint, the production chain of the sugarcane biomass facility consists of the following key elements:
\begin{enumerate}[label=\roman*.]
    \item \textbf{Products} $\rightarrow$ Includes inputs for the production process, such as sugarcane and straw, as well as intermediate and final outputs such as electricity, hydrogen, sugar, ethanol, and bagasse. We denote the complete set of products as $\mathcal{I}$;
    \item \textbf{Processes} $\rightarrow$ These are sequences of operations that convert input materials (e.g., sugarcane and straw) into intermediate/final products (e.g., electricity or biochar). The set of processes is referred to as $\mathcal{J}$;
    \item \textbf{Plants} $\rightarrow$ Production units that incorporate one or more processes. The set of plants is referred to as $\mathcal{K}$.
\end{enumerate}

\noindent To ease presentation, hereinafter, the set of all processes that can be executed in a production plant $k \in \mathcal{K}$ will be denoted by $\mathcal{J}(k)$, where $\cup_{k \in \mathcal{K}} \mathcal{J}(k) = \mathcal{J}$. Next, we formally introduce and discuss the proposed investment (first-stage) and operation (second-stage) models.

\subsection{Operation (Second-Stage) Model} \label{Sec::OperModel}

Let $\mathbf{c} \triangleq \big\{ c_{k} \big\}_{k \in \mathcal{K}}$ be the vector of planned production capacities for each of the $k \in \mathcal{K}$ plants in the biomass production facility and $\tilde{\boldsymbol{\xi}} \triangleq \big\{\tilde{P}_{i}, \tilde{D}^{(0)}_{i}\big\}_{i \in \mathcal{I}}$ the uncertain vector containing the sales price $\tilde{P}_{i}$ of the products $i \in \mathcal{I}^{(sell)} \subseteq \mathcal{I}$ which can be sold on the market and the initial availability $\tilde{D}^{(0)}_{i}$ of all products $i \in \mathcal{I}$ at the start of the second stage. Each process $j \in \mathcal{J}$ involves a set of input products $\mathcal{I}^{\text{(in)}}(j) \subset \mathcal{I}$, and a set of output products $\mathcal{I}^{\text{(out)}}(j) \subset \mathcal{I}$. Similarly, let $\mathcal{J}^ {\text{(in)}}(i) = \big\{ j \in \mathcal{J} ~ | ~ i \in \mathcal{I}^{\text{(in)}}(j)\big\}$ be the set of all processes for which product $i \in \mathcal{I}$ is an input and $\mathcal{J}^{\text{(out)}}(i) = \big\{ j \in \mathcal{J} ~ | ~ i \in \mathcal{I}^{\text{(out)}}(j)\big\}$ the set of all processes such that $i \in \mathcal{I}$ is an output product.

Firstly, the proportions of products involved in each process are given by known constants $\theta^{\text{(in)}}_{j, i}, ~ \forall ~ i \in \mathcal{I}^{in}(j)$ and $\theta^{\text{(out)}}_{j, i}, ~ \forall i ~ \in \mathcal{I}^{\text{(out)}}(j)$. Thus, to ensure that the quantities of each input $f^{\text{(in)}}_{j, i}, ~ \forall i \in \mathcal{I}^{\text{(in)}}(j)$ and output $f^{\text{(out)}}_{j, i}, \forall ~ i \in \mathcal{I}^{\text{(out)}}(j)$ products follow these proportions, an auxiliary variable $l_{j}$ is introduced, representing the ``production level'' of process $j \in \mathcal{J}$. Thus, for a given process $j \in \mathcal{J}$, the following relationship must hold:
\begin{align}
    & f^{\text{(in)}}_{j,i}  = l_{j} \theta^{\text{(in)}}_{j, i} \geq 0,  && \forall ~ i \in \mathcal{I}^{\text{(in)}}(j); \label{eq:2nd-c0} \\
    & f^{\text{(out)}}_{j,i} = l_{j} \theta^{\text{(out)}}_{j, i} \geq 0, && \forall ~ i \in \mathcal{I}^{\text{(out)}}(j). \label{eq:2nd-c1}
\end{align}

\noindent Additionally, the model maintains product balance across the production system. For each product $i \in \mathcal{I}$, the final quantity $d^{(f)}_{i}$ is determined by an initial (uncertain) availability $\tilde{D}^{(0)}_{i}$, the amount consumed $f^{\text{(in)}}_{j, i}$ in each process $j \in \mathcal{J}^{\text{(in)}}(i)$, the amount produced $f^{\text{(out)}}_{j,i}$ in each process $j \in\mathcal{J}^{\text{(out)}}(i)$, and the quantity sold in the market $v_{i} \geq 0$. The balance equation for each product $i \in \mathcal{I}$ is thus expressed as follows:
\begin{align}
    & v_{i} \le \tilde{D}^{(0)}_{i} + \sum_{j \in\mathcal{J}^{\text{(out)}}(i)} f^{\text{(out)}}_{j,i} - \sum_{j \in \mathcal{J}^{\text{(in)}}(i)} f^{\text{(in)}}_{j,i}. \label{eq:2nd-c2}
\end{align}

\noindent The operation of each plant $k \in \mathcal{K}$ is limited by its planned capacity $c_{k}$. Following operational standards, the capacity of a plant is measured by the total volume of \textit{reference products} processed across all processes that can be executed on it. More precisely, let $\phi(j) \in \mathcal{I}^{in}(j)$ denote the mapping that identifies the \textit{reference product} of a given process $j \in \mathcal{J}(k)$ of a production plant $k \in \mathcal{K}$. Therefore, the capacity constraint for a $k \in \mathcal{K}$ can be expressed as follows:
\begin{align}
    & \sum_{j \in \mathcal{J}(k)} f^{\text{(in)}}_{j, \phi(j)} \leq c_{k}, && \forall ~ k \in \mathcal{K} \label{eq:2nd-cap}.
\end{align}

Finally, the operation (second-stage) model aims to minimize the net cost, defined as the operational costs of each process offset by sales revenue. Specifically, the operational cost of a process $j \in \mathcal{J}$ is determined by the unit cost $O_{j}$ of its \textit{reference product} $\phi(j)$. On the other hand, sales revenue is calculated based on the uncertain sales price $\tilde{P}_{i}$ and volume sold $v_{i}$ of each product $i \in \mathcal{I}^{\text{(sell)}}$. Therefore, given a capacity plan $\mathbf{c}$ and uncertain vector $\tilde{\boldsymbol{\xi}}$, the operation model is defined as follows:
\begin{align}
    Q\big(\mathbf{c}, \tilde{\boldsymbol{\xi}}\big) = & \min_{\substack{\boldsymbol{f}^{\text{(in)}}, \boldsymbol{f}^{\text{(out)}}, \boldsymbol{d}^{(f)}, \\ \boldsymbol{v} \geq \boldsymbol{0}}} ~ \sum_{j \in \mathcal{J}} O_{j} f^{\text{(in)}}_{j, \phi(j)} - \sum_{i \in \mathcal{I}^{(sell)}} \tilde{P}_{i} v_{i} \label{eq01::2ndStageModel} \\
    & ~ \text{subject to: Constraints \eqref{eq:2nd-c0}--\eqref{eq:2nd-cap}} \label{eq02::2ndStageModel}.
\end{align}

For expository purposes and to ease reference throughout this work, Fig. \ref{fig:Illust_Prod} illustrates the production chain of an example facility composed of $\mathcal{K} = \big\{1, 2, 3\big\}$ plants and involving a total of $|\mathcal{I}| = 7$ products. For simplicity, each plant is assumed to contain a single process $\mathcal{J}(k) = \big\{k\big\}, ~ \forall ~ k \in \mathcal{K}$. Specifically, \textit{Plant} $k = 1$, called the \textit{Integrated Sugar Mill}, runs the \textit{Ethanol and Sugar Production} process. This process takes \textit{Sugarcane} (Product $i = 1$) as input and produces \textit{Sugar} (Product $i = 3$), \textit{Ethanol} (Product $i = 4$), and \textit{Bagasse} (Product $i = 5$) as outputs. Among these, only \textit{Sugar} and \textit{Ethanol} are final products sold on the market (represented by the triangle in the diagram). \textit{Bagasse} is an intermediate product used as input for the next plant (\textit{Plant} $k = 2$). Also, \textit{reference product} that characterizes the production cost and capacity of \textit{Plant} $k = 1$ is \textit{Sugarcane}, marked with a ${}^{*}$ in the diagram. \textit{Plant} $k = 2$ uses \textit{Bagasse} (from \textit{Plant} $k = 1$) and \textit{Straw} (product $i = 2$) as inputs in the \textit{Electricity from Residues} process, generating \textit{Electricity} (Product $i = 6$). \textit{Electricity} can be sold on the market and used as input for the next plant (\textit{Plant} $k = 3$). Furthermore, note that this process in \textit{Plant} $k = 2$ uses two input products: \textit{Bagasse}, a byproduct from \textit{Plant} $k = 1$, and \textit{Straw}, which serves as a raw feedstock. Given the characteristics of the \textit{Electricity from Residues} process in \textit{Plant} $k = 2$, \textit{Bagasse} is the \textit{reference product} for determining its production cost and capacity. Finally, \textit{Plant} $k = 3$ operates the \textit{Electrolysis} process, converting \textit{Electricity} into \textit{Hydrogen} (product $i = 7$), which is also sold on the market. For the \textit{Electrolysis} process, the \textit{reference product} is \textit{Electricity}. Leveraging the operation model \eqref{eq01::2ndStageModel}--\eqref{eq02::2ndStageModel}, in the next section, we present and discuss the proposed capacity expansion model.

\subsection{Capacity Expansion (First-Stage) Model} \label{Sec::ExpModel}

The primary goal of the investment model is to determine a capacity expansion plan for a sugarcane biomass production facility that balances investment and risk-adjusted operational costs. To ease presentation, let $y_{k} \geq 0$ be an auxiliary variable representing the expansion over existing capacity $C_{k}^{(0)}$ for a given plant $k \in \mathcal{K}$. Thus, the total capacity planned for plant $k \in \mathcal{K}$ is defined as 
\begin{align}
    c_{k} = C_{k}^{(0)} + y_{k} \geq 0. \label{eq:1st-c0}
\end{align}

From a financial standpoint, investment costs typically do not scale linearly with new capacity levels. In fact, this well-known phenomenon of ``economies of scale'' is particularly evident within the sugarcane industry. For instance, \cite{tribe} examined and validated the applicability of the capacity–cost scaling law originally proposed by \cite{chilton}, which describes how capital expenditures scale with plant capacity. According to this rule, if the capacity of a planned plant is $\beta$ times that of a reference plant, the required CAPEX increases by a factor of $\beta^{\sigma}$, where a typical value of $\sigma$ is 0.7. In this context, if $\beta = 10$, the total investment would be approximately $10^{0.7} \approx 5$ times that of the reference plant, meaning the cost per unit of output is halved. Conversely, if the planned capacity is $\beta = 0.1$ (i.e., ten times smaller), the investment would be around $0.1^{0.7} \approx 0.2$ times the reference cost, resulting in a doubling of unit cost (${0.2}/{0.1} = 2$). Naturally, if some production capacity already exists, the associated investment is subtracted from the total investment requirement. 

Therefore, by acknowledging this fact, in this work, we assume that the total capital cost $\big(b_{k}\big)$ of the production plant $k \in \mathcal{K}$ is nonlinearly associated with new $(y_{k})$ and existing $\big(C_{k}^{(0)}\big)$ capacities. Formally, following the discussion in \cite{tribe, chilton}, let $B^{(\text{ref})}_{k}$ be the \textit{a priori} known reference cost of a given candidate for expansion plant $k \in \mathcal{K}$ with reference capacity $C^{(\text{ref})}_{k}$, and $\sigma_{k}$ the aforementioned scaling exponent associated with economies of scale. Then, for each candidate plant $k \in \mathcal{K}$, the relationship between total investment cost and capacity is given by the following expression.
\begin{align}
    & b_{k} = B_{k}^{(\text{ref})} \left( \frac{C_{k}^{(0)} + y_{k}}{C_k^{(\text{ref})}} \right)^{\sigma_{k}} - B_{k}^{(\text{ref})} \left( \frac{C_{k}^{(0)}}{C_{k}^{(\text{ref})}} \right)^{\sigma_{k}} \geq 0 \label{eq:1st-c1}
\end{align}

In this work, the investment plan adopts the standard time-based approach, in which capacity expansion is planned to operate during a representative (typical) year \cite{Tanaka2022_EconomicsofPwrSyst}. Consequently, for consistency in the analysis, both capital and operational costs must be expressed on an annual basis. In line with economic literature, this is achieved using the concept of an annuity. To convert the total investment expenditure $b_{k}$ into an annualized cost, let $N_{k}$ denote the useful lifetime (in years) of plant $k \in \mathcal{K}$ and $R_{k}$ the opportunity cost of capital for investing in its expansion. The annuity coefficient is then calculated as follows:
\begin{align}
    & \varphi_{k} = \frac{R_{k}}{1 - (1 + R_{k})^{-N_{k}}}, && \forall ~ k \in \mathcal{K}.
\end{align}

Accordingly, to capture the investor's risk aversion, the proposed capacity expansion model employs the $\alpha$-quantile-based risk measure CVaR \cite{Rockafellar2002_CVaRforGeneralLossDist, Barbosa2024_TailoredWIndOp, Street2025_ManagingRESRisk_Brasil_HourlyHedge} and is formulated as follows:
\begin{align}
    & \hspace{-0.00cm} \min_{\substack{\boldsymbol{b}, \mathbf{c}, \\ \mathbf{y} \geq \boldsymbol{0}}} ~ \sum_{k \in \mathcal{K}} \varphi_{k} b_{k} + \lambda ~ CVaR_{\alpha}\Big(Q(\mathbf{c}, \tilde{\boldsymbol{\xi}})\Big) + (1-\lambda) \mathbb{E}\Big[Q(\mathbf{c}, \tilde{\boldsymbol{\xi}})\Big] \label{eq01::1stStageModel} \\
    & \hspace{-0.00cm} \> \text{subject to: Constraints \eqref{eq:1st-c0}--\eqref{eq:1st-c1}}. \label{eq02::1stStageModel}
\end{align}

The solution to model \eqref{eq01::1stStageModel}-\eqref{eq02::1stStageModel} provides the risk-adjusted optimal expansion plan for the sugarcane facility under uncertain prices and product availability. However, model \eqref{eq01::1stStageModel}-\eqref{eq02::1stStageModel} is a nonlinear optimization problem, which is often difficult for standard solvers to handle efficiently. To address this computational challenge, the next section details a reformulation procedure. This procedure uses a piecewise linear approximation to create a sample-based version of the model that can be solved efficiently with existing optimization tools.

\section{Solution Method: Piecewise Approximation and Scenario-Based Representation} \label{sec:adress-problems}

The capacity expansion planning model presented in Section \ref{Sec::CapExpModel} follows a two-stage decision-making framework, where the first stage is framed as a nonlinear optimization problem due to the economies of scale relationship between total capacity and CAPEX defined in \eqref{eq:1st-c1}. To enable a computationally efficient formulation to be solved using standard optimization algorithms, two key challenges must be addressed: (i) handling the nonlinear relationship in \eqref{eq:1st-c1}, and (ii) devising a tractable representation of uncertain factors. Section \ref{sec:linearization} introduces an efficient piecewise linear approximation of the economies of scale effect using binary variables \cite{Dantzig1960_SigLinearProgProb_Integer, Lin2013_ReviewPiecewiseLinMethods, vielma} and Section \ref{sec:samples} then reformulates problem \eqref{eq01::1stStageModel}–\eqref{eq02::1stStageModel} as an equivalent scenario-based model.

\subsection{Piecewise Linear Approximation of Economies of Scale Effect} \label{sec:linearization}

As discussed in Section \ref{Sec::ExpModel}, capital expenditure typically does not scale linearly with capacity. In the context of the sugarcane industry, this relationship can be effectively approximated using the capacity–cost scaling law \cite{tribe, chilton}, which is mathematically represented in \eqref{eq:1st-c1}. While this nonlinear formulation adds realism to the capacity expansion planning model, it also introduces computational challenges. To address this issue, we adopt the piecewise linear approximation approach proposed by \cite{Dantzig1960_SigLinearProgProb_Integer}, which uses a binary-variable-based representation to approximate nonlinear relationships. Formally, for each plant $k \in \mathcal{K}$, let $\big\{(c_{k,p}, b_{k,p})\big\}_{p = 1}^{M+1}$ be a predefined, ordered set of $(M + 1)$ capacity–cost pairs that follow the nonlinear relationship specified in \eqref{eq:1st-c0}–\eqref{eq:1st-c1}. Then, the nonlinear relationship between total capacity and investment expenditure over the interval $\big[c_{k,1}, c_{k,M+1}\big] \times \big[b_{k,1}, b_{k,M+1}\big]$ can be approximated by the following mixed-integer set of piecewise linear functions \cite{Dantzig1960_SigLinearProgProb_Integer}:
\begin{align}
    & \hspace{-0.60cm} \Pi\bigg( \Big\{(c_{k,p}, b_{k,p})\Big\}_{p = 1}^{M+1} \bigg) = \Bigg\{ \big(c_{k}, b_{k}\big) \in \mathbb{R}^{2}_{+} ~ \Bigg| \notag \\
    & \hspace{2.00cm} \exists ~ \big(\boldsymbol{\gamma}, \boldsymbol{\eta}\big) \in \big[0, 1\big]^{M} \times \big\{0, 1\big\}^{M}, ~~ \text{s.t.} \notag \\
    & \hspace{2.00cm} c_{k} = \sum_{p = 1}^{M} \eta_{p} \Big(\gamma_{p} c_{k,p} + \big(1 - \gamma_{p}\big) c_{k,p+1}\Big); \notag \\
    & \hspace{2.00cm} b_{k} = \sum_{p = 1}^{M} \eta_{p} \Big(\gamma_{p} b_{k,p} + \big(1 - \gamma_{p}\big) b_{k,p+1}\Big); \notag \\
    & \hspace{2.00cm} \sum_{p = 1}^{M} \eta_{p} \leq 1 \Bigg\}. \label{eq::PISet}
\end{align}

\begin{figure}[!ht]
    \centering
    \includegraphics[width=\linewidth]{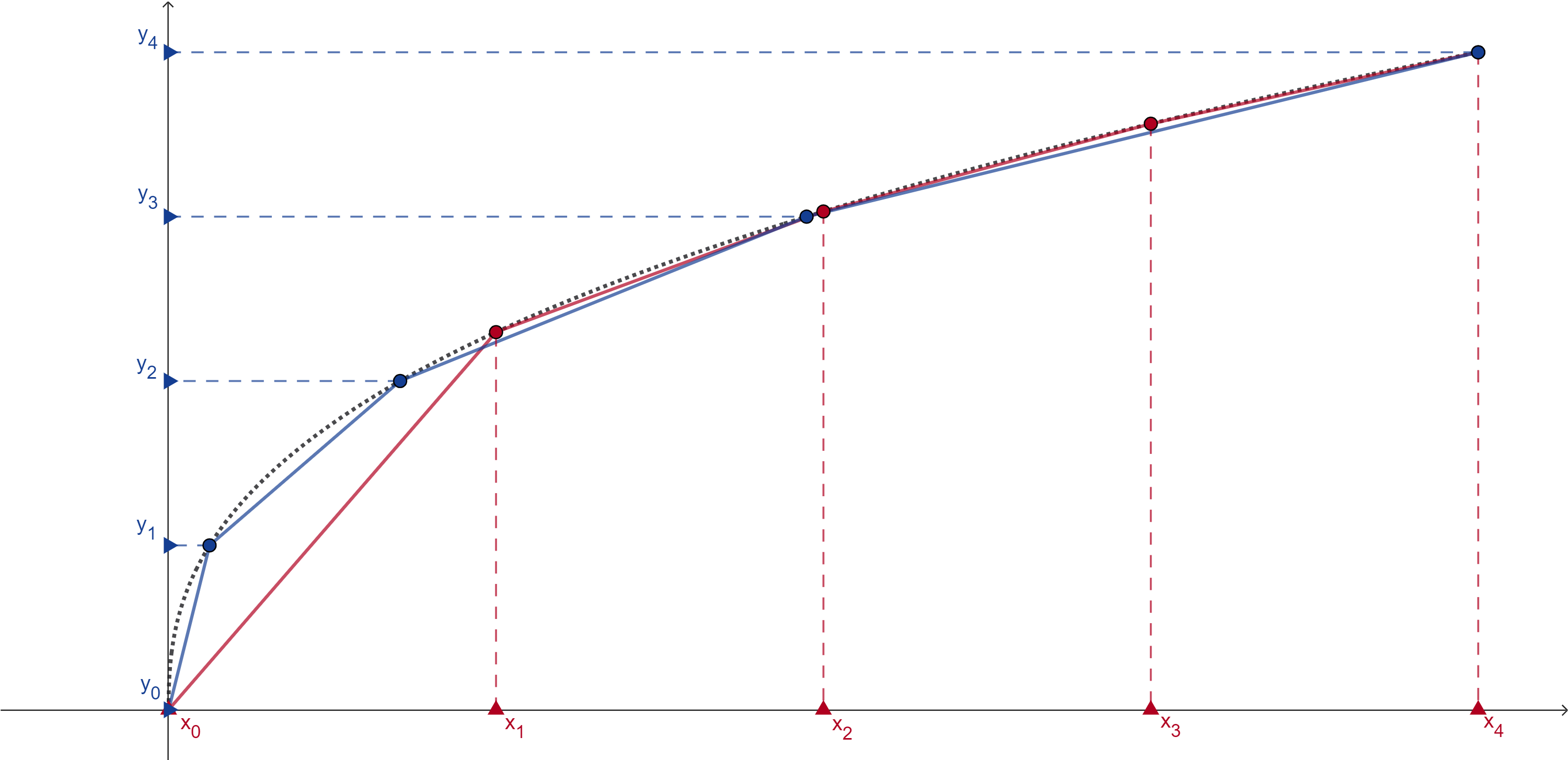}
    \caption{Piecewise linear approximation with points uniformly distributed along the \myred{domain} and the \myblue{range}.}
    \label{fig:linearizations}
\end{figure}

It is important to note from \eqref{eq::PISet} that the accuracy of the piecewise linear approximation increases with the number of segments $M$. However, this improvement comes at the cost of greater computational complexity, as evaluating the set $\Pi$ becomes more demanding for higher values of $M$. Therefore, a balance between accuracy and computational efficiency must be considered by the decision-maker. In the particular context of this work, the resulting capital cost functional is concave in total capacity. Thus, to enhance approximation accuracy without introducing an excessive number of binary variables, we adopt a discretization strategy that uniformly partitions the ``range'' space $\big[b_{k,1}, b_{k,M+1}\big]$ and computes the corresponding values in the ``domain'' space $\big[c_{k,1}, c_{k,M+1}\big]$ using \eqref{eq:1st-c0}–\eqref{eq:1st-c1}. This approach results in a denser distribution of points where the cost-capacity curve changes rapidly (i.e., at lower capacities) and a sparser distribution where the curve becomes nearly linear (i.e., at higher capacities). Fig. \ref{fig:linearizations} illustrates the difference in approximation quality when points are uniformly distributed across the domain versus the range. By leveraging the mixed-integer set \eqref{eq::PISet} to approximate the nonlinear relationship between total capacity and investment cost, this approach captures the economies of scale effect more accurately while remaining computationally tractable, enabling the use of state-of-the-art MILP solvers, which have seen significant advances in performance over recent decades \cite{Koch2022_ProgressMPSolvers_2021_2020}.

\subsection{Equivalent Scenario-Based Capacity Expansion Model} \label{sec:samples} 

To address the two-stage structure of problem \eqref{eq01::1stStageModel}–\eqref{eq02::1stStageModel}, we adopt the standard approach widely used in both the scientific literature and industry to represent the uncertainty through a finite set of scenarios $\big(\Omega\big)$. Specifically, the uncertain factors $\big\{\tilde{P}_{i}, \tilde{D}^{(0)}_{i}\big\}_{i \in \mathcal{I}}$ are described by a sample set $\big\{\big(\{P_{i,\omega}, D^{(0)}_{i,\omega}\}_{i \in \mathcal{I}}, \pi_{\omega}\big)\big\}_{\omega \in \Omega}$, where $\boldsymbol{\pi} \triangleq \big\{\pi_{1}, \dots, \pi_{|\Omega|}\big\}$ denotes the associated probability distribution over scenarios. In this context, since the second-stage problem is a linear program, and leveraging the CVaR formulation from \cite{Rockafellar2002_CVaRforGeneralLossDist}, the original two-stage problem \eqref{eq01::1stStageModel}–\eqref{eq02::1stStageModel} can be equivalently reformulated as the following scenario-based mathematical programming problem:
\begin{align}
    & \min_{\substack{\boldsymbol{b}, \mathbf{c}, \mathbf{y}, \boldsymbol{\gamma}, \boldsymbol{\eta}, \\ \boldsymbol{f}^{\text{(in)}}_{\omega}, \boldsymbol{f}^{\text{(out)}}_{\omega}, \boldsymbol{v}_{\omega}}} ~ \sum_{k \in \mathcal{K}} \varphi_{k} b_{k} + \lambda \Bigg(z + \frac{1}{1 - \alpha} \sum_{\omega \in \Omega} \pi_{\omega} \delta_{\omega}\Bigg) \notag \\
    & \hspace{4.60cm} + \big(1 - \lambda\big) \sum_{\omega \in \Omega} \pi_{\omega}Q_{\omega} \label{eq01::SingLevel} \\
    & ~ \text{subject to:} \notag \\
    & ~ f^{\text{(in)}}_{j, i, \omega} = l_{j} \theta^{\text{(in)}}_{j, i} \geq 0 ,  && \hspace{-4.57cm} \forall ~ i \in \mathcal{I}^{\text{(in)}}(j),  ~ j \in \mathcal{J}, ~ \omega \in \Omega; \label{eq02::SingLevel} \\
    & ~ f^{\text{(out)}}_{j, i, \omega} = l_{j} \theta^{\text{(out)}}_{j, i} \geq 0, && \hspace{-4.70cm} \forall ~ i \in \mathcal{I}^{\text{(out)}}(j), ~ j \in \mathcal{J}, ~ \omega \in \Omega; \label{eq03::SingLevel} \\
    & ~ v_{i,\omega} \leq D^{(0)}_{i,\omega} + \sum_{j \in\mathcal{J}^{\text{(out)}}(i)} f^{\text{(out)}}_{j,i,\omega} - \sum_{j \in \mathcal{J}^{\text{(in)}}(i)} f^{\text{(in)}}_{j,i,\omega}, \notag \\
    & && \hspace{-2.60cm} \forall ~ i \in \mathcal{I}, ~ \omega \in \Omega; \label{eq04::SingLevel} \\
    & ~ \sum_{j \in \mathcal{J}(k)} f^{\text{(in)}}_{j, \phi(j), \omega} \leq c_{k}, && \hspace{-2.73cm} \forall ~ k \in \mathcal{K}, ~ \omega \in \Omega; \label{eq05::SingLevel} \\
    & ~ Q_{\omega} = \sum_{j \in \mathcal{J}} O_{j} f^{\text{(in)}}_{j, \phi(j), \omega} - \sum_{i \in \mathcal{I}^{(sell)}} P_{i,\omega} v_{i, \omega}, && \hspace{-1.55cm} \forall ~ \omega \in \Omega; \label{eq06::SingLevel} \\
    & ~ c_{k} = C_{k}^{(0)} + y_{k} \geq 0, && \hspace{-1.55cm} \forall ~ k \in \mathcal{K}; \label{eq07::SingLevel} \\
    & ~ \big(c_{k}, b_{k}\big) \in \Pi\bigg( \Big\{(c_{k,p}, b_{k,p})\Big\}_{p = 1}^{M+1} \bigg), && \hspace{-1.55cm} \forall ~ k \in \mathcal{K}; \label{eq08::SingLevel} \\
    & ~ z + \delta_{\omega} \geq Q_{\omega}, && \hspace{-1.55cm} \forall ~ \omega \in \Omega; \label{eq09::SingLevel} \\
    & ~ \delta_{\omega}, b_{k}, y_{k} \geq 0. && \hspace{-1.55cm} \forall ~ k \in \mathcal{K}. \label{eq11::SingLevel}
\end{align}

\noindent Problem \eqref{eq01::SingLevel}–\eqref{eq11::SingLevel} represents an implementable reformulation of the original two-stage program \eqref{eq01::1stStageModel}-\eqref{eq02::1stStageModel}, and belongs to the class of two-stage stochastic MILP models with recourse. In the numerical experiment presented in Section \ref{Sec::NumericalExperiment}, we solve this extensive form directly, without applying any decomposition techniques.

\section{Numerical Experiment} \label{Sec::NumericalExperiment}

To illustrate the applicability and effectiveness of the proposed capacity expansion planning model, a set of numerical experiments is conducted using realistic data from the Brazilian energy system. Specifically, Section \ref{sec:experimental-setup} outlines the experimental setup and base case for analysis. Section \ref{sec:case1} examines the investment plan under a risk-neutral perspective, while Section \ref{sec:case2} extends the analysis to a risk-averse profile. Finally, Section \ref{sec:case3} explores production routes that were found to be less advantageous in the base cases. Each execution of the model took less than 5 minutes, demonstrating the efficiency of the proposed approach.

\subsection{Experimental Setup} \label{sec:experimental-setup}

The base case assumes a producer manages 46,000 hectares of sugarcane with a productivity of 84 tons per hectare (the Brazilian average), yielding an annual capacity of 3.86 million tons. This integrated facility currently produces both sugar and 1G ethanol. To account for price uncertainty, we use a pure data-driven (nonparametric) approach. Price scenarios are constructed from observed monthly closing prices for international sugar (New York) and 1G ethanol (Chicago) indices, spanning from January 2008 to August 2024 (i.e., $|\Omega| = 200$ scenarios in total), with a probability density of $\pi_{\omega} = 1/|\Omega|, ~ \forall ~ \omega \in \Omega$ assigned to each observation. For the other products considered, sale prices were chosen to reflect recent conditions in the Brazilian market and were kept constant across scenarios.

The analysis explores new investment opportunities for the producer, including the expansion of its integrated sugar and ethanol plant and the potential development of additional facilities for other sugarcane-derived products. Beyond sugar and 1G ethanol, the portfolio of final products encompasses bioelectricity and advanced biofuels such as biomethane, hydrogen, bio-oil, SAF, green diesel, bio-gasoline, and methanol. The scope also includes biochar and carbon credits. In addition, the production pipeline incorporates by-products such as straw, bagasse, vinasse, filter cake, and molasses, which serve as inputs for other production routes. The $\text{CO}_{2}$ released during processing is treated as a by-product, with the potential to be converted into carbon credits.

For the integrated sugar and ethanol plant, two main processes are considered: one focus on sugar production, yielding 87 kilograms of sugar and 28 liters of ethanol per ton of sugarcane, and another focus on ethanol production, producing 43 kilograms of sugar and 57 liters of ethanol per ton of sugarcane. Both processes generate by-products such as bagasse, vinasse, and filter cake. Vinasse and filter cake can be used in biomethane production. Additionally, there is the possibility of investing in electricity generation from residues, using bagasse and straw as inputs. The analysis considers three processes for producing the trio of fuels: SAF, bio-gasoline, and green diesel. Two routes involve the AtJ plant: AtJ 1G and AtJ 2G, which are the same physical process but use either 1G or 2G ethanol as input. The third route is the Fischer-Tropsch process, which utilizes bagasse as a feedstock. Hydrogen can be produced using Steam Ethanol Reform (with both 1G and 2G ethanol) or Steam Methane Reform (which requires biomethane). A third option is electrolysis, using residue-generated electricity. Biochar and bio-oil are co-produced via pyrolysis in a single plant: slow pyrolysis favors biochar, while fast pyrolysis favors bio-oil. The facility can also earn carbon credits through BECCS (Bioenergy with Carbon Capture and Storage). Capturing the CO2 from the integrated sugar and ethanol process transitions production from carbon-neutral to carbon-negative, actively reducing atmospheric carbon. \answer{The technical details and data for this case study are described in Appendix \ref{sec:appendix} and in the external repository \cite{dataset_}. In addition, the repository also provides supplementary analytical results with visualization tools.}

\subsection{Case Study 1: Risk Neutral}\label{sec:case1}

In this case study, we assume a risk-neutral profile drives the investment plan. The resulting investment strategy involves adding an electricity generation plant that uses straw and bagasse to electricity production. The added capacity is sufficient to process part of the available straw (228 thousand of 541 thousand tons) and the entire bagasse output from sugar and ethanol production (966 thousand tons), enabling annual generation of 439 gigawatt-hours. More than half of the straw will be left in the field, contributing to soil fertility restoration, a practice traditionally adopted in sugarcane farming. The energy conversion is efficient using high-pressure boilers, optimizing steam generation and electricity production. Sugar and ethanol sales prices have a significant impact on plant operations. When the price of sugar per ton exceeds 65\% of the price of ethanol per kiloliter, the process prioritizes sugar production; otherwise, the ethanol-focused process is chosen. Sugar is preferred in 70.1\% of the scenarios. In scenarios where both commodity prices are low, the plant does not operate at full capacity, resulting in some sugarcane remaining unprocessed. Electricity production, in contrast, consistently utilizes all available bagasse, yielding an average of 421 gigawatt-hours per year. The process also generates vinasse and filter cake as by-products. Across scenarios, the average net revenue is 52.13 MM\$, with values ranging from –19.91 MM\$ to 180.79 MM\$. Notably, 16.5\% of scenarios result in negative net revenues, representing a significant financial risk for the producer.

\subsection{Case Study 2: Risk Averse}\label{sec:case2}

Case 2 builds on the base setup (Section \ref{sec:experimental-setup}) but introduces a risk-averse profile ($\lambda=0.5$ and $\alpha=0.9$) to guide investment. In addition to the new Electricity from Residues plant, an AtJ plant is added. This plant uses 1G ethanol to produce SAF, green diesel, and bio-gasoline, with a capacity sized to consume half of the producer's maximum ethanol output.
With these expansions, the facility processes all available sugarcane in every scenario, ensuring the maximum amount of bagasse is available and allowing electricity generation to reach its full capacity (439 gigawatt-hours per year) across the board.
The production of sugar versus ethanol remains sensitive to their price ratio, though sugar is prioritized in only 22.5\% of scenarios. Ethanol production is split, with the AtJ plant consuming a portion and the remainder being sold. The AtJ plant operates at full capacity in 74\% of the scenarios. However, when ethanol prices or the sugar-to-ethanol price ratio are high, the production of SAF, green diesel, and bio-gasoline is reduced or even stopped.
Consistent with the CVaR risk policy, the average net revenue remains high (51.00 MM\$), but the financial risk is significantly reduced: the range of outcomes is much narrower (-2.91 MM\$ to 167.71 MM\$), and the proportion of scenarios resulting in revenue losses falls sharply to nearly 1\%.

\subsection{Case Study 3: Economic Use of Waste}\label{sec:case3}

In the previous analysis, several products were not included in the optimal expansion and operation plan. Therefore, to further explore these products, this case study conducts a sensitivity analysis to identify the minimum sales prices at which production of selected products becomes economically feasible. For expository purposes, we use the base case data while adopting the risk-averse profile described in Section \ref{sec:case2}.

\subsubsection{Biomethane}\label{sec:biomethane}

In the base case analysis, the reference sales price of biomethane was set equal to the price of a natural gas molecule, which is currently 0.40 US dollars per cubic meter in Brazil. For economic viability, however, the sale price of biomethane would need to rise to at least 0.52 US dollars per cubic meter. At this price level, the biomethane plant processes all of the filter cake generated in sugar and ethanol production, along with part of the vinasse, yielding 16.02 MMm\textsuperscript{3} annually. It is worth noting that this price is close to the final tariff paid by industry. Nevertheless, this tariff reflects not only the cost of the molecule but also transport fees and the distributor's margin. When restricted to the molecule component alone, biomethane would need to be roughly 30\% more expensive to achieve viability.

\subsubsection{Hydrogen}

In the previous analyses, none of the three technological routes to produce hydrogen proved viable, as they relied on the price of fossil hydrogen, set at 1.50 US dollars per kilogram. For economic viability, the minimum sales price of hydrogen must reach 2.31 US dollars per kilogram. Although this level represents a premium over fossil hydrogen, it remains competitive with the water electrolysis (green hydrogen) route. More specifically, the optimal capacity expansion plan selects steam methane reforming as the pathway, with biomethane serving as the required input. Biomethane production capacity is the same as defined in Section \ref{sec:biomethane}, while the steam methane reforming facility has sufficient capacity to process the full biomethane output, yielding 3.8 thousand tons of hydrogen annually. Both facilities operate at full capacity across all scenarios.

\subsubsection{Methanol}

As a more easily transportable fuel, has diverse applications, including use as a marine fuel, biodiesel precursor, and feedstock in the chemical industry. The technological routes considered in these numerical experiments involve methanol synthesis using either hydrogen or biomethane as inputs. In the base case analysis, the assumed sales price of methanol reflects the value of fossil methanol, set at 400 \$/ton. To achieve economic viability, the sales price would need to rise to 1045 \$/ton, approximately three times higher than the base case value.

\subsubsection{Biochar}

Biochar is a particularly valuable product since, when applied to land, it can boost crop productivity by up to 15\%. As the model does not capture cycles within the production chain, the sensitivity analysis for biochar was carried out by fixing a biomass pyrolysis plant with sufficient capacity to produce the required biochar for land application. The harvesting process was then adjusted to reflect the productivity gains from biochar use. As a result, sugar and ethanol production capacity expanded to process the additional sugarcane, while electricity generation capacity also increased due to the higher availability of bagasse. Revenues were adjusted by deducting the annual cost of the pyrolysis plant (introduced with an initial capacity) and by excluding potential revenues from biochar sales, since the product was assumed to be applied to land rather than sold in the market. After these adjustments, the average net revenue rose to 69.35 MM\$, representing a 36\% improvement over the base case. This gain, which exceeds the direct increase in crop yield, can be attributed to economies of scale.

\section{Conclusion}

    This work proposes a two-stage stochastic optimization model to guide capacity expansion planning for sugarcane biomass facilities, integrating long-term investment and short-term operational decisions under uncertainty in market prices and resource availability. \answer{The model explicitly accounts for diverse production pathways, including electricity, hydrogen, sugar, ethanol, advanced biofuels, biomethane, and biochar. On the one hand, it captures economies of scale through a capacity–cost scaling relationship, leading to a nonconvex optimization problem that is addressed through an efficient and tractable reformulation as an MILP model. This treatment constitutes a key differentiating feature relative to standard approaches, which typically simplify capacity expansion modeling by assuming constant marginal investment costs. On the other hand, the framework incorporates the investor's risk preferences through a CVaR-based functional that balances expected revenue and downside risk, thereby enabling a robust trade-off analysis among profitability, diversification, and market volatility. Overall, the strength of the proposed framework lies in its ability to jointly optimize multiple dimensions of the supply chain, allowing emerging and waste-to-energy technologies to be evaluated alongside traditional products within a unified risk-averse decision-making structure. In addition to its methodological contributions, this work also provides an open-source implementation of the proposed capacity expansion planning framework, referred to as \textit{OptBio} \cite{OptBio}. By making the model and solution workflow publicly available, we promote transparency, reproducibility, and practical adoption, allowing researchers and practitioners to readily apply, extend, and benchmark the proposed approach.} A numerical case study based on the Brazilian energy system demonstrates the model's strategic utility. Results show that risk-neutral strategies favor traditional products but lead to high volatility. Conversely, risk-averse strategies prioritize diversification, particularly with investments in Alcohol-to-Jet processes, significantly improving the risk-return profile. Finally, sensitivity analyses indicate that waste-derived products (biomethane, hydrogen) can become viable under specific prices, and biochar enhances crop productivity, indirectly boosting electricity and ethanol output.

    \answer{As a direction for future research, we highlight exploring a robust-optimization-based formulation of the capacity expansion planning model developed in this work, with the objective of mitigating the impact of uncharacterizable variations in product prices on the resulting investment decisions. Furthermore, characterizing the uncertainty of both product prices and feedstock availability through a joint, multivariate modeling framework is highly valuable to both the academic community and practitioners. In addition, developing a computationally efficient characterization of concave economies-of-scale cost functions that further enhances the computational performance of the proposed capacity expansion model is a worthwhile direction for future research. Extending the proposed capacity expansion framework to a broader context that explicitly incorporates environmental objectives represents a promising direction for future research. In particular, adopting a central planner's perspective to analyze investment equilibria aimed at achieving predefined environmental targets through alternative investment–operation strategies constitutes a natural extension of the framework developed in this work. Finally, future work could explore the integration of demand-side flexibility, following approaches such as \cite{Du2025_Decarb_DataCenterNet_PowerMitigation}; incorporate endogenous sources of uncertainty, as in \cite{Giannelos2025_StochOptiModel_Expans_ExoEndoUncert, Pianco2025_DDU_DRO_ExpansionDistSystem}; and scale the modeling and analysis to regional or national levels by coordinating the investment and operation of multiple heterogeneous facilities, in line with the ideas proposed in \cite{Dong2024_FlexEnhanc_UbanSystem}.}

\section*{Declaration of generative AI and AI-assisted technologies in the writing process}

During the preparation of this work the author(s) used GPT-4 to enhance the clarity of select portions of the text. After using these tools, the author(s) reviewed and edited the content as needed and take(s) full responsibility for the content of the published article.

\bibliographystyle{IEEEtran}
\bibliography{bibliography.bib}

\onecolumn

\newpage

\appendices

\section{} \label{sec:appendix}

\begin{landscape}    
        \begin{figure}
            \centering
            \includegraphics[width=1\linewidth]{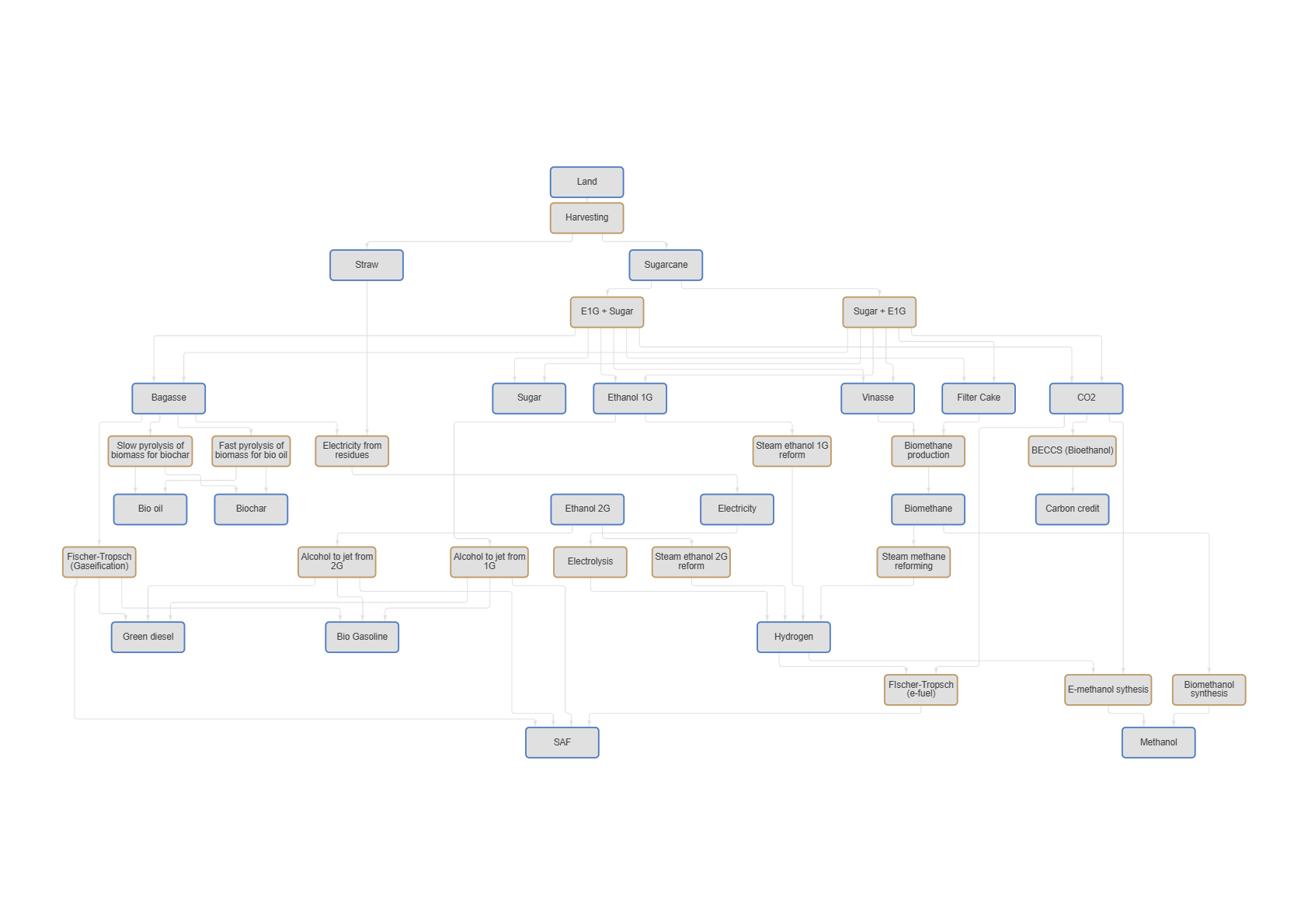}
            \caption{Schematic representation of the production chain considered in the case studies.}
            \label{fig:ProcessDiagram}
        \end{figure}
    \end{landscape}

    \answer{This appendix provides a detailed description of the experimental setup underlying the numerical analysis presented in Section \ref{Sec::NumericalExperiment}. Table \ref{tab:plants} summarizes the plants considered in the study, including their existing (initial) capacities, as well as the corresponding reference CAPEX and capacity values. Table \ref{tab:products} presents the full spectrum of products considered in the experiment, including production inputs as well as intermediate and final outputs, along with their key characteristics, such as initial availability and selling price. Table \ref{tab:processes} provides a detailed description of the processes considered in the case study, highlighting their OPEX values, input–output information, and the corresponding plant to which each process is associated. Finally, Fig. \ref{fig:ProcessDiagram} presents a schematic representation of the production chain, showing the processes and their relationships with the products. Boxes outlined in blue represent products, whereas boxes outlined in orange represent processes. Arrows indicate the flow of products through the processes. To facilitate data access and reproducibility, the external repository \cite{dataset_} contains these data and technical details, as well as supplementary analytical results and visualization tools.}
    
    \begin{center}
        \begin{table}[!ht]
        \centering
        \caption{Summary of the plants considered in the numerical experiment in Section \ref{Sec::NumericalExperiment}, including their existing (initial) capacities, as well as the corresponding reference CAPEX and capacity values.}
        \renewcommand{\arraystretch}{1.30}
        \label{tab:plants}
        \small
        \setlength{\tabcolsep}{6pt}
            \begin{tabular}{l l c c c}
                \hline
                \textbf{Plant} & \textbf{Capacity unit} & \textbf{Initial Capacity} & \textbf{Reference CAPEX (MM\$)} & \textbf{Reference Capacity} \\
                \hline
                Harvesting & kha of Land & 46.0 & 0.0 & 500.0 \\
                Biomethane production & kt of Filter Cake & 0.0 & 48.0 & 115.8 \\
                Steam methane reforming & Mm$^{3}$ of Biomethane & 0.0 & 3.5 & 58.8 \\
                Pyrolysis of biomass & kt of Bagasse & 0.0 & 52.0 & 84.0 \\
                Alcohol-to-jet & ML of Ethanol 1G & 0.0 & 97.7 & 109.5 \\
                Fischer--Tropsch (Gasification) & kt of Bagasse & 0.0 & 1330.0 & 966.0 \\
                Electricity from residues & kt of Straw & 0.0 & 29.9 & 250.0 \\
                Ethanol 1G + Sugar mill & Mt of Sugarcane & 3.87 & 196.6 & 4.38 \\
                Electrolysis & GWh of Electricity & 0.0 & 249.0 & 832.2 \\
                BECCS (Bioethanol) & Mt of CO$_2$ & 0.0 & 1009.0 & 15.7 \\
                E-methanol synthesis & kt of CO$_2$ & 0.0 & 34.3 & 110.0 \\
                Fischer--Tropsch (e-fuel) & kt of Hydrogen & 0.0 & 77.0 & 26.7 \\
                Biomethanol production & Mm$^{3}$ of Biomethane & 0.0 & 79.0 & 18.8 \\
                Steam ethanol 1G reforming & ML of Ethanol 1G & 0.0 & 38.1 & 123.4 \\
                Steam ethanol 2G reforming & ML of Ethanol 2G & 0.0 & 38.1 & 123.4 \\
                \hline
            \end{tabular}
        \end{table}
    \end{center}

    \begin{table}[!ht]
        \centering
        \caption{Full spectrum of products considered in the experiment, including production inputs as well as intermediate and final outputs, along with their key characteristics, such as initial availability and selling price.}
        \label{tab:products}
        \small
        \setlength{\tabcolsep}{6pt}
        \renewcommand{\arraystretch}{1.30}
        \begin{tabular}{l l r r}
            \hline
            \textbf{Product} & \textbf{Unit} & \textbf{Initial Availability} & \textbf{Average Sell Price} \\
            \hline
            Land & ha & 46000.0 & 25.25 \\
            Sugarcane & t & 0.0 & 19.49 \\
            Straw & t & 0.0 & 25.25 \\
            Sugar & t & 0.0 & 389.74 \\
            Molasses & t & 0.0 & 194.89 \\
            Filter Cake & t & 0.0 & 0.00 \\
            Bagasse & t & 0.0 & 0.00 \\
            Ethanol 1G & L & 0.0 & 0.51 \\
            Ethanol 2G & L & 0.0 & 1.02 \\
            Vinasse & t & 0.0 & 0.00 \\
            Electricity & MWh & 0.0 & 40.0 \\
            Biomethane & m$^{3}$ & 0.0 & 0.40 \\
            Hydrogen & t & 0.0 & 1500.00 \\
            Bio oil & t & 0.0 & 500.00 \\
            SAF & L & 0.0 & 1.85 \\
            Green diesel & L & 0.0 & 1.70 \\
            Bio Gasoline & L & 0.0 & 1.60 \\
            Methanol & t & 0.0 & 400.00 \\
            Biochar & t & 0.0 & 1000.00 \\
            CO$_2$ & t & 0.0 & 0.00 \\
            Carbon credit & t & 0.0 & 30.00 \\
            \hline
        \end{tabular}
        \end{table}

    \begin{table}[!ht]
        \centering
        \caption{Detailed description of the processes considered, highlighting their OPEX values, input–output information, and the corresponding plant to which each process is associated.}
        \label{tab:processes}
        \small
        \setlength{\tabcolsep}{6pt}
        \renewcommand{\arraystretch}{1.30}
        \begin{tabular}{l r l l l}
            \hline
            \textbf{Process} & \textbf{Opex} & \textbf{Input} & \textbf{Output} & \textbf{Plant} \\
            \hline
            Harvesting & 1200.0 & Land: 1.0 ha & Sugarcane: 84.0 t & Harvesting \\
             &  &  & Straw: 11.76 t &  \\
            \hline
            Biomethane production & 16.0 & Filter Cake: 115800.0 t & Biomethane: 1.6e7 m$^{3}$ & Biomethane production \\
             &  & Vinasse: 1.35872e6 t &  &  \\
            \hline
            Steam ethanol 1G reform & 0.01 & Ethanol 1G: 8.29 L & Hydrogen: 0.001 t & Steam ethanol 1G Reforming \\
            \hline
            Steam ethanol 2G reform & 0.01 & Ethanol 2G: 8.29 L & Hydrogen: 0.001 t & Steam ethanol 2G Reforming \\
            \hline
            Steam methane reforming & 0.02 & Biomethane: 4.2 m$^{3}$ & Hydrogen: 0.001 t & Steam methane reforming \\
            \hline
            Fast pyrolysis of biomass for bio oil & 73.3 & Bagasse: 1.0 t & Bio oil: 0.55 t & Pyrolysis of biomass \\
             &  &  & Biochar: 0.2 t &  \\
            \hline
            Slow pyrolysis of biomass for biochar & 49.0 & Bagasse: 1.0 t & Bio oil: 0.352 t & Pyrolysis of biomass \\
             &  &  & Biochar: 0.47 t &  \\
            \hline
            Alcohol-to-jet from 1G & 0.21 & Ethanol 1G: 109.5 L & SAF: 40.3 L & Alcohol-to-jet \\
             &  &  & Green diesel: 9.8 L &  \\
             &  &  & Bio Gasoline: 5.7 L &  \\
            \hline
            Alcohol-to-jet from 2G & 0.21 & Ethanol 2G: 109.5 L & SAF: 40.3 L & Alcohol-to-jet \\
             &  &  & Green diesel: 9.8 L &  \\
             &  &  & Bio Gasoline: 5.7 L &  \\
            \hline
            Fischer--Tropsch (Gasification) & 204.7 & Bagasse: 1.0 t & SAF: 90.68 L & Fischer--Tropsch (Gasification) \\
             &  &  & Green diesel: 46.07 L &  \\
             &  &  & Bio Gasoline: 36.13 L &  \\
            \hline
            Biomethanol synthesis & 0.17 & Biomethane: 817.4 m$^{3}$ & Methanol: 1.0 t & Biomethanol production \\
            \hline
            Electricity from residues & 12.28 & Straw: 70.0 t & Electricity: 135.0 MWh & Electricity from residues \\
             &  & Bagasse: 297.0 t &  &  \\
            \hline
            Electrolysis & 5.0 & Electricity: 60.0 MWh & Hydrogen: 1.0 t & Electrolysis \\
            \hline
            BECCS (Bioethanol) & 2.57 & CO$_2$: 1.0 t & Carbon credit: 1.0 t & BECCS (Bioethanol) \\
            \hline
            E-methanol synthesis & 578.38 & CO$_2$: 1.26 t & Methanol: 1.0 t & E-methanol syntesis \\
             &  & Hydrogen: 146.0 t &  &  \\
            \hline
            Fischer--Tropsch (e-fuel) & 0.6 & Hydrogen: 1.0 t & SAF: 0.75 L & Fischer--Tropsch (e-fuel) \\
             &  & CO$_2$: 1.0 t &  &  \\
            \hline
            Sugar + E1G & 24.0 & Sugarcane: 1.0 t & Sugar: 0.0867 t & Ethanol 1G + Sugar Mill \\
             &  &  & Bagasse: 0.25 t &  \\
             &  &  & Filter Cake: 0.03 t &  \\
             &  &  & CO$_2$: 0.0133 t &  \\
             &  &  & Ethanol 1G: 28.333 L &  \\
             &  &  & Vinasse: 310.66 t &  \\
            \hline
            E1G + Sugar & 24.0 & Sugarcane: 1.0 t & Sugar: 0.04333 t & Ethanol 1G + Sugar Mill \\
             &  &  & Bagasse: 0.25 t &  \\
             &  &  & Vinasse: 621.33 t &  \\
             &  &  & Filter Cake: 0.03 t &  \\
             &  &  & CO$_2$: 0.026667 t &  \\
             &  &  & Ethanol 1G: 56.67 L &  \\
            \hline
        \end{tabular}
    \end{table}

\end{document}